# 13 или 14, Хэллоуиновский сюрприз
Алексей Панов, Алексей Хованский

Ромбокубооктаэдр – это один из высокосимметричных Архимедовых многогранников. Лаконичное описание Архимедовых многогранников[1], их еще называют Архимедовыми телами, содержится в трактате Паппа Александрийского *Собрание* (*Συναγωγή*), опубликованном где-то около 340 года. В 1588 году он был переведен с греческого на латынь и с тех пор стал доступен европейским математикам.

Иоганн Кеплер был знаком с этим переводом. В своем великом произведении *Гармония Мир* (*Harmonices Mundi*), над которым он работал с 1599 по 1619 год, Кеплер дает развернутое определение Архимедовых тел. Используя его, он доказывает существование ровно тринадцати Архимедовых тел, присваивает каждому из них собственное название и впервые предъявляет полный набор всех их изображений[2].

**Ромбокубооктаэдр**. Некоторые из этих тринадцати многогранников были известны математикам и художникам Возрождения, предшественникам Кеплера. Вот, например, две иллюстрации Леонардо да Винчи 1509 года, на которых ромбокубооктаэдр изображен в двух разных стилях (рис. 1).

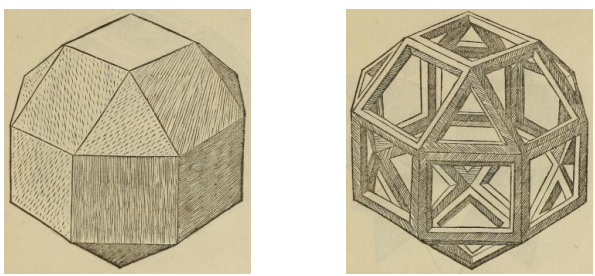

Рис. 1. Ромбокубооктаэдр, Леонардо да Винчи 1509

Можно сосчитать, что у него 26 граней, из которых 8 – правильные треугольники и 18 – квадраты.

Займёмся сборкой модели ромбододекаэдра. Возьмем лист плотной бумаги, например лист ватмана формата А2, и вырежем из него три заготовки. Одна из них предназначена для экваториального пояса модели – это полоска из 9 одинаковых квадратов со стороной 5 см (рис. 2).

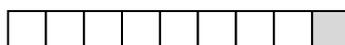

Рис. 2. Заготовка для экваториального пояса

Две другие – для полярных областей модели, каждая из таких же 9 квадратов (рис. 3).

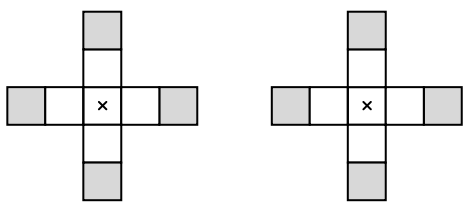

Рис. 3. Заготовки для полярных областей

---

[1] [Chris Rorres, Archimedean Solids](#)
[2] [Johannes Kepler, Harmonice Mundi, pp. 61-65](#)

На обоих рисунках 2, 3 серым цветом выделены поверхности квадратов, предназначенные для склейки, и еще на рисунке 3 крестиками указаны положения будущих полюсов.

Теперь на каждой заготовке нужно сделать сгибы по сторонам соседних квадратов, склеить полоску в восьмигранный экваториальный пояс и расположить все три заготовки как на рисунке 4.

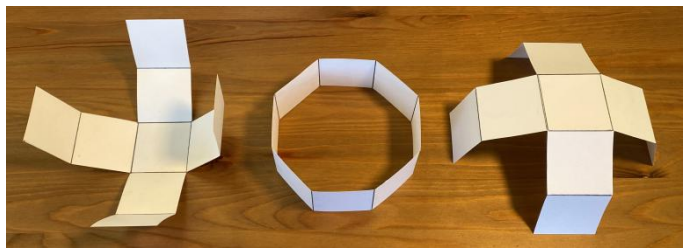

Рис. 4. Три заготовки для ромбокубооктаэра

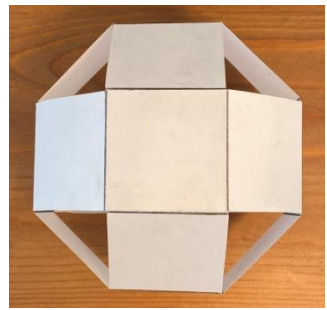

Рис. 5. Ромбокубооктаэдр

А потом полярные области нужно сверху и снизу подклеить к экваториальному поясу. Окончательный результат – на рисунке 5. Это гибридная модель ромбокубооктаэдра – нечто среднее между левой и правой частью рисунка 1. Можно сказать, что треугольные грани у этой модели прозрачные.

По тому, как мы собирали ромбокубооктаэдр, можно подумать, что у него один экваториальный пояс и пара полюсов. Но теперь, когда у нас есть готовая модель, легко убедиться, что на самом деле у него 3 экваториальных пояса и, соответственно, 3 пары полюсов.

Используя эту модель, можно заметить, что при вращении вокруг некоторых прямых на некоторые углы ромбокубоокаэдр совмещается сам с собой. Эти прямые называются осями поворотной симметрии ромбокубооктаэдра. Оказывается, что их у него ровно 13, это все прямые, соединяющие центры его противоположных граней. Например, ромбокубооктаэдр самосовмещается при поворотах на 90°, 180° и 270° вокруг полярной оси, то есть вокруг прямой, проходящей через два его противоположных полюса. В любом случае, можно сказать, что ромбокубооктаэдр – это высокосимметричный многогранник.

**Псевдоромбокубооктаэдр, № 14?** Псевдоромбокубооктаэдр – это младший родственник ромбокубооктаэдра, он появился на математическом горизонте только в прошлом веке, но у него имеется своя предыстория, о которой мы расскажем чуть позже. А сейчас у нас уже есть все необходимое для его сборки. Вернемся к рисунку 4 и повернем одну из полярных областей на нем на 45° (рис. 6).

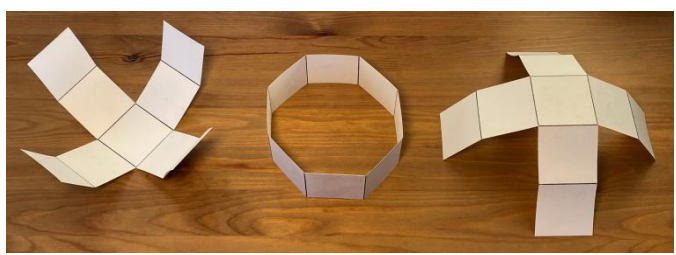

Рис. 6. Три заготовки для псевдоромбокубооктаэра

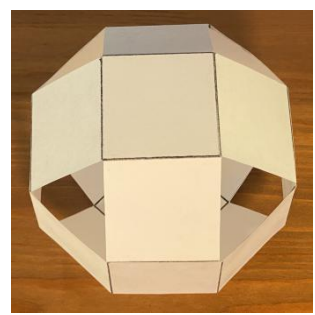

Рис. 7. Псевдоромбокубооктаэдр

После этого точно так же, как раньше, подклеим сверху и снизу полярные области к экваториальному поясу, и псевдоромбокубооктаэдр готов (рис. 7)

Полученный многогранник похож на ромбокубооктаэдр, у него те же самые 8 треугольных и 18 квадратных граней, тоже имеется экваториальный пояс с полярными областями. Но только он намного менее симметричен. У него, например, всего лишь 5 осей поворотной симметрии.

В своем эссе *О шестиугольных снежинках*, изданном в 1611 году, во время работы над *Гармонией мира*, Кеплер пишет

> *Вспомнив о ромбах, я приступил к геометрическим изысканиям, чтобы выяснить, какое тело, аналогичное пяти правильным и **четырнадцати** Архимедовым телам, можно составить из одних ромбов.[3]*

Упоминание 14 Архимедовых тел во время работы над *Гармонией мира* и над своим полным списком 13 Архимедовых тел ни в коем случае не может быть опиской Кеплера. Возможно, к моменту написания эссе Кеплер как раз обнаружил псевдоромбокубооктаэдр и посчитал его еще за одно Архимедово тело. И только чуть позже, признав псевдоромбокубооктаэдр недостаточно симметричным, изменил свою точку зрения.

**Хэллоуинский сюрприз**. А теперь поменяем дизайн нашей модели псевдоромбокубооктаэдра. Склеим ее из плотной оранжевой бумаги и снабдим двумя красными светодиодами с батарейками[4].

И вот окончательный результат.

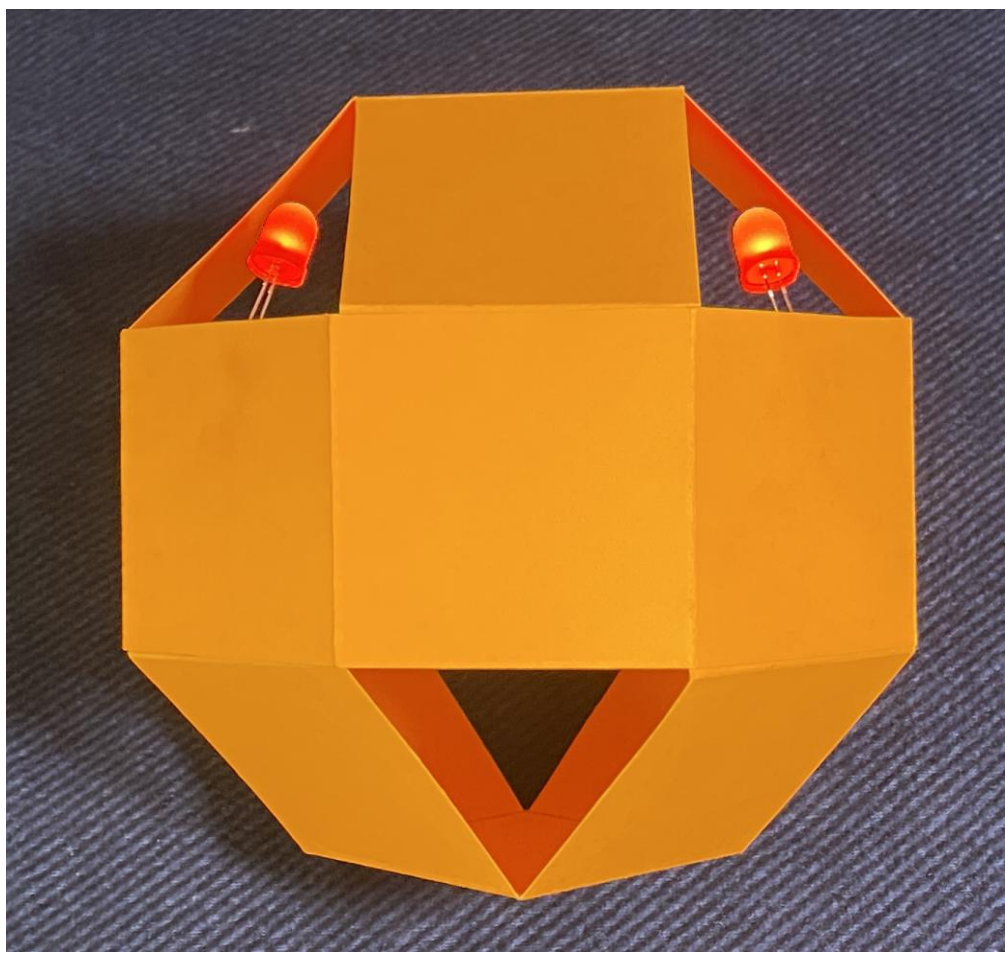

Рис. 8. Хэллуиновский сюрприз

---

[3] Иоганн Кеплер, О шестиугольных снежинках. — М.: Наука, 1982, с. 10
[4] Дмитрий Панов, Фонарик-светлячок. Квантик 2014, № 12